\documentclass[11pt]{article}
\usepackage{amssymb,amsmath,latexsym, bbm}

\usepackage{epsfig}

 \newtheorem{thm}{Theorem}[section]

\newtheorem{cor}[thm]{Corollary}

\newtheorem{lemma}[thm]{Lemma}

\newcommand{\pf}{\noindent{\bf Proof.} }
\def\qed{{\hfill $\Box$ \bigskip}}

\def\R{\mathbb{R}}

\def\Z{\mathbb{Z}}

\def\P{\mathbb{P}}

\def\E{\mathbb{E}}

\def\p{\partial}

\def\<{\langle} \def\>{\rangle}

\def\mK{{\mathbb K}}

\def\LL{\mathcal{L}}
\def\sL{\mathcal{L}}
\def\SS{\mathcal{S}}

\def\1{\mathbbm{1}}

\def\wh{\widehat}
\def\wt{\widetilde}
\def\eps{\varepsilon}

\def\bee{\begin{equation}}
\def\eee{\end{equation}}

\makeatletter
\@addtoreset{equation}{section}

\makeatother

\begin{document}
\bibliographystyle{plain}

\title{\Large \bf
Heat Kernels for Non-symmetric Non-local Operators
}

\author{{\bf Zhen-Qing Chen}
\thanks{Research partially supported by NSF Grant
 DMS-1206276.}
\quad and \quad {\bf Xicheng Zhang}
\thanks{Research partially supported by NNSFC grant of China (Nos. 11271294, 11325105).}
}
 
\date{}
\maketitle

\begin{abstract} We survey the recent progress in the study of
heat kernels  for a class of non-symmetric non-local operators.
  We focus on the existence and
sharp two-sided estimates  of the heat kernels and their connection
to jump diffusions.
\end{abstract}

\bigskip
\noindent {\bf AMS 2000 Mathematics Subject Classification}: Primary
60J35, 47G30, 60J45; Secondary: 31C05, 31C25, 60J75

\bigskip

 \noindent{\bf Keywords}:  Discontinuous Markov process, diffusion with jumps, non-local operator, pseudo-differential
operator, heat kernel estimate,  L\'evy system

\section{Introduction}

Second order elliptic differential operators and diffusion processes
take up, respectively, an central place in the theory of partial
differential equations (PDE) and the theory of probability. There
are close relationships between these two subjects. For a large
class of second order elliptic differential operators $\LL$ on
$\R^d$, there is a diffusion process $X$ on $\R^d$ associated with
it so that $\LL$ is the infinitesimal generator of $X$, and vice
versa.  The connection between $\LL$ and $X$ can also be seen as
follows. The fundamental solution (also called heat kernel)
  for $\LL$ is the transition density function of $X$.
For example, when
$$
\LL=\frac12 \sum_{i, j=1}^d a_{ij}(x) \frac{\partial^2}{\partial x_i \partial x_j} + \sum_{i=1}^d b_i (x) \frac{\partial}{\partial x_i},
$$
 where   $(a_{ij}(x))_{1\leq i, j\leq d}$  is a $d\times d$ symmetric
 matrix-valued  continuous function on $\R^d$ that is uniformly elliptic and bounded,
and $b (x) =(b_1 (x), \dots , b_d(x))$ is  a bounded $\R^d$-valued function,
there is a unique diffusion
$X=\{X_t, t\geq 0; \P_x, x\in \R^d\}$
on $\R^d$ that solves the martingale problem for $(\LL, C_c^2 (\R^d))$.
That is, for every $x\in \R^d$, there is a unique probability measure $\P_x$  on the space $C([0, \infty); \R^d)$
of continuous functions on $\R^d$ so that 
$\P_x(X_0=x)=1$ and for every $f\in C^2_c (\R^d )$,
$$
M^f_t:= f (X_t)-f(X_0)- \int_0^t \LL f (X_s) ds
$$
 is a $\P_x$-martingale.
Here $X_t (\omega )=\omega (t)$ is the coordinate map on $C([0, \infty); \R^d)$.
It is also known that $(X, \P_x)$ is the unique weak solution to
the following stochastic differential equation
$$
dX_t= \sigma (X_t) dW_t + b(X_t) dt, \qquad X_0=x,
$$
where $W_t$ is an $n$-dimensional Brownian motion and 
$\sigma (x)= a(x)^{1/2}$ is the symmetric square root matrix of  $a(x)=(a_{ij}(x))_{1\leq i,j\leq d}$.

 When $a$ is 
H\"older continuous,
it is known that  $\LL$ has a jointly continuous heat kernel $p(t, x, y)$
with respect to the Lebesgue measure on $\R^d$ that enjoys
the following Aronson's estimate (see Theorem \ref{T:4.3} below):
there are constants $c_k>0$,
 $k=1, \cdots, 4$, so that
\begin{equation}\label{e:1.1}
 c_1 t^{-d/2}\exp(-c_2|x-y|^2/t) \leq p(t, x, y)\leq  c_3 t^{-d/2}\exp(-c_4|x-y|^2/t)
 \end{equation}
 for $t>0$ and $x, y \in \R^d$.

As many physical and economic
systems exhibit discontinuity or jumps,
in-depth study on
non-Gaussian jump  processes are called for.
See for example, \cite{B, JW, KSZ, STa} and the references therein.
 The infinitesimal generator of a discontinuous Markov process in $\R^d$
is no longer a differential operator but rather a non-local (or,
integro-differential) operator. For instance, the infinitesimal
generator of an isotropically
 symmetric $\alpha$-stable process in $\R^d$ with
$\alpha \in (0, 2)$ is a fractional Laplacian operator $c\,
\Delta^{\alpha /2}:=- c\, (-\Delta)^{\alpha /2}$. During the past
  several years  there
  is also many interest
  from the theory of PDE (such as singular obstacle problems) to study
non-local operators; see, for example,
 \cite{CSS, S} and the references
therein. Quite many progress has been made in the last fifteen years on
the  development of the
   De Giorgi-Nash-Moser-Aronson type  theory for   non-local operators.
   For example,     Kolokoltsov  \cite{kol} obtained two-sided heat kernel
estimates for certain stable-like processes in $\R^d$, whose
infinitesimal generators are   a class of pseudo-differential
operators having smooth symbols.
 Bass and Levin
\cite{BL02b} used a completely different approach to obtain similar
estimates for discrete time Markov chain on $\Z^d$, where the
conductance between $x$ and $y$ is comparable to $|x-y|^{-n-\alpha}$
for $\alpha \in (0,2)$. In \cite{CK1}, two-sided heat kernel
estimates and a scale-invariant parabolic Harnack inequality (PHI in abbreviation) for
symmetric $\alpha$-stable-like processes on $d$-sets are obtained.
 Recently in \cite{CK2},
two-sided heat kernel estimates and PHI are
established for symmetric non-local operators of variable order.
The De Giorgi-Nash-Moser-Aronson type  theory is studied very recently in
\cite{CK3} for symmetric diffusions with jumps.
  We refer the reader to the survey articles \cite{C, GHL} and the references therein on the study of 
  heat kernels for symmetric non-local operators.
   However,   for non-symmetric non-local operators,  much less is known.
   In this article, we will survey the recent development in the study of heat kernels
   for non-symmetric non-local operators. We will concentrate on the recent progress
   made in  \cite{CZ, CZ2} and \cite{CHXZ}.
   In Section \ref{S:5} of this paper, we summarize
some other recent work on heat kernels for non-symmetric non-local operators.
We also take this opportunity to fill a gap in the proof of \cite[(3.20)]{CZ},
which is \eqref{e:3.23a} of this paper. The proof in \cite{CZ} works for the case 
$|x|\geq t^{1/\alpha}$. In Section 3, a proof is supplied for the case $|x|\leq t^{1/\alpha}$.
 In fact,  a slight modification of the original proof for \cite[Theorem 2.5]{CZ} 
gives  a better estimate   \eqref{e:3.18} than \eqref{e:3.23a}.
 
  In this survey, we concentrate on heat kernel on the whole Euclidean spaces and on
  the work that the authors are involved.
  We will not discuss Dirichlet heat kernels in this article.

\section{L\'evy process}

A L\'evy process on $\R^d$ is a right continuous process
$X=\{X_t; t\geq 0\}$ having left limit
that has independent stationary increments. It is uniquely characterized
by its L\'evy exponenent $\psi$:
\begin{equation}\label{e:2.1}
\E_0 \exp (i\xi \cdot X_t) = \exp ( -t \psi (\xi)), \quad \xi \in \R^d.
\end{equation}
Here for $x\in \R^d$,
the subscript $x$ in the mathematical expectation $\E_x$ and the probability $\P_x$ means that the process $X_t$ starts from $x$.
The L\'evy exponent $\psi$ admits a unique decomposition:
\begin{equation}\label{e:2.2}
\psi (\xi )= i b\cdot \xi+  \sum_{i, j=1}^d a_{ij} \xi_i\xi_j
+ \int_{\R^d}
\left( 1-e^{i\xi \cdot z} + i\xi \cdot z \1_{\{ |z|\leq 1\}} \right)
\Pi (dz) ,
\end{equation}
where $b\in \R^d$ is a constant vector, $(a_{ij})$ is a non-negative definite symmetric constant matrix, and  $\Pi (dz)$
is a positive measure
on $\R^d \setminus \{0\}$ so that
$\int_{\R^d} (1\wedge |z|^2) \Pi (dz) <\infty$.
The L\'evy measure $\Pi (dz)$ has a strong probabilistic meaning.
It describes the jumping intensity of $X$ making a jump of size $z$.
 Denote by
 $\{P_t; t\geq 0\}$ the transition semigroup of $X$;
 that is, $P_t f (x) = \E_x f(X_t)
=\E_0 f(x+X_t)$.
For an integrable function $f$,
its Fourier transform is defined to be $\wh f(\xi )=\int_{\R^d} e^{i\xi \cdot x} f(x) dx$. Then we have by \eqref{e:2.1} and Fubini's theorem,
$$
\wh {P_t f} (\xi) =\int_{\R^d} e^{i\xi \cdot x} \E_0 f(x+X_t) dx
= \E_0 \left[ e^{-\xi \cdot X_t} \left( \int_{\R^d} e^{i\xi \cdot (x+X_t)}
f(x+X_t) dx \right) \right]
= e^{-t \psi (-\xi)} \wh f (\xi).
$$
If we denote the infinitesimal generator of $\{P_t; t\geq 0\}$ (or $X$)
by $\LL$,
then
\begin{equation}\label{e:2.3}
\wh {\LL f}(\xi) =\frac{d}{dt}\Big|_{t=0} \wh {P_t f}(\xi) = -\psi (-\xi) \wh f (\xi).
\end{equation}
Hence $-\psi (-\xi)$ is the Fourier multiplier (or symbol) for the infinitesimal generator $\LL$ of $X$. One can derive a more explicit expression for the generator $\LL$: for $f\in C^2_c(\R^d)$,
\begin{equation}\label{e:2.4}
\LL f(x) =\sum_{i,j=1}^d a_{ij} \frac{\partial^2 f }{\partial x_i
\partial x_j} (x)+ b\cdot \nabla f(x)
+ \int_{\R^d} \left( f(x+z)-f(x) -\nabla f(x) \cdot z \1_{\{|z|
\leq 1\}} \right) \Pi (dz).
\end{equation}

When $b=0$, $\Pi=0$ and 
 $(a_{ij})={\bf I}_{d\times d}$ 
the identity matrix,
that is when $\psi (\xi)= |\xi|^2$, $X$ is a Brownian motion in $\R^d$
with variance $2t$ and infinitesimal generator $\Delta:=\sum_{i=1}^d
\frac{\partial^2}{\partial x_i^2}$. When $b=0$, $a_{ij}=0$ for all $1\leq i, j\leq d$ and
 $\Pi (dz)= 
 {\cal A}(d, -\alpha)
  |z|^{-(d+\alpha)}d z$
for $0<\alpha <2$,
where ${\cal A}(d, -\alpha)$
 is a normalizing constant so that
$\psi (\xi )= |\xi|^\alpha$, $X$ is a rotationally symmetric
$\alpha$-stable process in $\R^d$, whose infinitesimal generator
is the fractional Laplacian $\Delta^{\alpha/2}:= -(-\Delta)^{\alpha/2}$.

Unlike Brownian motion case, explicit formula for the transition density function
of symmetric $\alpha$-stable processes is not known except for a very few cases.
However we can get its two-sided estimates as follows.
It follows from \eqref{e:2.1} that under $\P_0$,
(i) $AX_t$ has the same distribution as $X_t$ for every $t>0$ and   rotation $A$ (an orthogonal matrix);
(ii) for every $\lambda >0$, $X_{\lambda t}$ has the same distribution
as $\lambda^{1/\alpha} X_t$. Let $p(t, x)$ be the density function of $X_t$ under $\P_0$; that is,
$$
p(t, x)= (2\pi)^{-d} \int_{\R^d} e^{-i x\cdot \xi} e^{-t|\xi|^\alpha} d \xi.
$$
Then $p(t, x)$ is a function of $t$ and $|x|$ and $p(t, x)=t^{-d/\alpha} p(1, t^{-1/\alpha} x)$.   Using Fourier's inversion,
one gets
$$
\lim_{ |x| \to \infty} |x|^{d+\alpha} p(1, x)= \alpha 2^{\alpha-1} \pi^{-(d/2+1)} \sin ( {\alpha \pi}/2 )
\Gamma ((d+\alpha)/2) \Gamma ( {\alpha}/2 ).
$$
(See P\'olya \cite{Po} when $d=1$ and Blumenthal-Getoor \cite[Theorem 2.1]{BG} when $d\geq 2$.)
It follows that $p(1, x) \asymp 1 \wedge  \frac{1}{|x|^{d+\alpha}} $.  Consequently,
\begin{equation}\label{e:2.5}
p(t, x) \asymp
  t^{-d/\alpha} \wedge  \frac{t}{|x|^{d+\alpha}}  \asymp \frac{t}{(t^{1/\alpha} + |x|)^{d+\alpha}}.
\end{equation}
Here for $a, b\in \R$, $a\wedge b:=\min\{a, b\}$,  and for two functions $f, g$, $f\asymp g$ means that
$f/g$ is bounded between two positive constants.

 In real world, 
almost every media we 
encounter
has impurities so we need to consider state-dependent 
 stochastic processes 
and state-dependent local and non-local operators. Intuitively speaking, we need to consider processes and operators where $\psi (\xi)$  is dependent on $x$;
that is, $\psi (x, \xi)$. If one uses Fourier multiplier approach \eqref{e:2.3}, one gets pseudo differential operators. 
The connection between pseudo differential operators and Markov processes has been nicely exposited in N. Jacob \cite{Ja}. 
 In this survey, 
we take \eqref{e:2.4} as a starting
point but with $a_{ij}(x)$, $b(x)$ and $\Pi (x, dz)$ being functions of $x\in \R^d$. That is,
$$
\LL f(x) =\sum_{i,j=1}^d a_{ij}(x) \frac{\partial^2 f }{\partial x_i
\partial x_j} (x)+ b(x) \cdot \nabla f(x)
+ \int_{\R^d} \left( f(x+z)-f(x) -\nabla f(x) \cdot z \1_{\{|z|
\leq 1\}} \right) \Pi (x, dz).
$$

We will concentrate on the case where $\Pi (x, dz)= \frac{\kappa (x, z)}{|z|^{d+\alpha}}dz$
for some $\alpha \in (0, 2)$ and  a measurable function $\kappa (x, z)$ 
on $\R^d\times \R^d$ satisfying for any $x, y, z\in \R^d$,
\begin{equation}\label{e:2.6}
0<\kappa_0\leq k(x, z)\leq \kappa_1 <\infty,
\qquad \kappa (x, z)=\kappa (x, -z),
\end{equation}
 and for some $\beta \in (0, 1)$,
\begin{equation}\label{e:2.7}
| \kappa(x, z) - \kappa (y, z) | \leq \kappa_2 |x-y|^\beta.
\end{equation}

\section{Stable-like processes and their heat kernels}

In this section, we consider the case where $a_{ij}=0$, $b=0$
and $\Pi (x, dz)=\frac{\kappa (x, z)}{|z|^{d+\alpha}}dz$; that is,
\begin{equation}\label{e:3.1}
\LL f(x) ={\rm p.v.} \int_{\R^d} \left( f(x+z)-f(x)  \right)  \frac{\kappa (x, z)}{|z|^{d+\alpha}}  dz,
\end{equation}
where $\kappa(x, z)$ is a function on $\R^d\times \R^d$ satisfying   \eqref{e:2.6} and \eqref{e:2.7}.
 Here p.v. stands for the Cauchy principal value, that is,
$$
\LL f(x) =\lim_{\eps \to 0} \int_{\{z\in \R^d: |x| \geq \eps\}}
\left( f(x+z)-f(x)  \right)  \frac{\kappa (x, z)}{|z|^{d+\alpha}}  dz.
$$
Since $\kappa (x, z)$ is symmetric in $z$, when $f\in C^2_b(\R^d)$, we can rewrite $\LL f (x)$ as
\begin{equation}\label{e:3.2}
\LL f(x) =\int_{\R^d} \left( f(x+z)-f(x) -\nabla f(x) \cdot z \1_{\{|z|
\leq 1\}} \right)  \frac{\kappa (x, z)}{|z|^{d+\alpha}}  dz.
\end{equation}
The non-local operator $\LL$ of \eqref{e:3.1} typically is not symmetric, 
as oppose to non-local operator given by
\begin{equation}\label{e:3.3}
\wt \LL f(x):= 
\lim_{\eps \to 0} \int_{\{y\in \R^d: |y-x| \geq \eps\}}
\left( f(y)-f(x)  \right)  \frac{c (x, y)}{|x-y|^{d+\alpha}}  dz 
\end{equation}
in the distributional sense. Here $c(x, y)$ is a symmetric function that is bounded between
two positive constants. The operator 
$\wt \LL$  
is the infinitesimal generator of the symmetric
$\alpha$-stable-like process studied in Chen and Kumagai \cite{CK1}, where it is shown that 
$\wt \LL$  
has a  jointly H\"older continuous 
heat kernel that admits  two-sided estimates
in the same form as \eqref{e:2.5}. 

The following result is recently established in \cite{CZ}.

\begin{thm} {\rm (\cite[Theorem 1.1]{CZ})} \label{T:3.1} 
Under \eqref{e:2.6} and \eqref{e:2.7}, 
there exists a unique nonnegative jointly continuous function $p  (t,x,y)$ in $(t,x,y)\in (0,1]\times\R^d\times\R^d$
solving
\begin{align}
\frac{\p }{\p t} p (t,x,y)=\LL p (t,\cdot,y)(x), \ x\not=y,  \label{eq15}
\end{align}
and enjoying
the following four properties:
\begin{description}
\item{\rm (i)}  (Upper bound)
There is a constant $c_1>0$
so that for all $t\in(0,1]$ and $x,y\in\R^d$,
\begin{align}
p (t,x,y)\leq c_1 t ( t^{1/\alpha}+|x-y|)^{-d-\alpha}. \label{eq16}
\end{align}

\item{\rm (ii)}  (H\"older's estimate)
For every $\gamma \in(0,\alpha\wedge 1)$, there is a constant
$c_2>0$ so that
for every
$t\in(0,1]$ and $x, ,y, z\in\R^d$,
\begin{align}
|p (t,x,z)-p (t, y, z)|\leq c_2|x-y|^\gamma t^{1-({\gamma}/{\alpha})}\left(t^{1/\alpha}+|x-z|\wedge|y-z|
\right)^{-d-\alpha} . \label{Ho}
\end{align}

\item{\rm (iii)}  (Fractional derivative estimate)
For all $x\not=y\in\R^d$,
the mapping $t\mapsto \LL p(t,\cdot,y)(x)$ is continuous on $(0,1]$, and
\begin{align}
|\LL p (t,\cdot,y)(x)|\leq c_3(t^{1/\alpha}+|x-y|)^{-d-\alpha}.\label{eq17}
\end{align}
\
\item{\rm (iv)}  (Continuity) For any bounded and uniformly continuous function $f:\R^d\to\R$,
\begin{align}\label{eq14}
\lim_{t\downarrow 0}\sup_{x\in\R^d}\left|\int_{\R^d}p (t,x,y)f(y)d y-f(x)\right|=0.
\end{align}
\end{description}
Moreover,  we have the following conclusions.

\begin{description}
\item{\rm (v)}  The constants $c_1$, $c_2$ and $c_3$ in {\rm (i)-(iii)} above
can be chosen so
that they  depend only
on $(d, \alpha, \beta, \kappa_0, \kappa_1, \kappa_2)$,  $(d, \alpha, \beta, \gamma, \kappa_0, \kappa_1, \kappa_2)$, and $(d, \alpha, \beta, \kappa_0, \kappa_1, \kappa_2)$, respectively.

\item{\rm (vi)}  (Conservativeness)
For all $(t,x,y)\in(0,1]\times\R^d\times\R^d$,
$p(t, x, y)\geq 0$ and
\begin{align}
\int_{\R^d}p (t,x,y)d y=1.\label{EK1}
\end{align}

\item{\rm (vii)}   (C-K equation)  For all $s,t\in(0,1]$ with $t+s\in(0,1]$
and $x,y\in\R^d$,
the following Chapman-Kolmogorov equation holds:
\begin{align}\label{eq201}
\int_{\R^d}p (t,x,z)p (s,z,y)d z=p (t+s,x,y).
\end{align}

\item{\rm (viii)}   (Lower bound)
 There exists $c_4=c_4(d, \alpha, \beta, \kappa_0, \kappa_1, \kappa_2)>0 $ so that for
 all $t\in(0,1]$ and $x,y\in\R^d$,
\begin{align}
p (t,x,y)\geq c_4t(t^{1/\alpha}+|x-y|)^{-d-\alpha}. \label{e:3.10}
\end{align}

\item{\rm (ix)}  (Gradient estimate)
For $\alpha\in[1,2)$,  there exists $c_5=c_5(d, \alpha, \beta, \kappa_0, \kappa_1, \kappa_2)>0 $
so that for all
 $x\not= y$ in $\R^d$
 and $t\in(0,1]$,
\begin{align}
|\nabla_x \log  p (t,x ,y)|\leq c_5t^{-1/\alpha}  .
\label{eq170}
\end{align}
\end{description}
\end{thm}

\subsection{Approach}

We now sketch the main idea behind the proof of Theorem \ref{T:3.1}.

To emphasize the dependence of $\LL$ in \eqref{e:3.1} on $\kappa$, we write it
as $\LL^\kappa$.
For each fixed $y\in \R^d$, we consider L\'evy process (starting from 0) with
 L\'evy measure $\Pi_y (dz)=\frac{\kappa (y, z)}{|z|^{d+\alpha}} dz$,
and denote its marginal probability density function and  infinitesimal generator by
 $p_y(t, x)$ and $\LL^{\kappa (y)}$, respectively.
Then we have
\begin{equation} \label{e:3.12}
\frac{\p }{\p t} p_y (t, x)= \LL^{\kappa (y)} p_y (t, x).
\end{equation}
We use Levi's idea and search for heat kernel  $p(t, x, y)$ for $\LL^\kappa $ with the following form:
\begin{equation}\label{e:3.13}
p(t, x, y)= p_y (t, x-y)+ \int_0^t \int_{\R^d} p_z (t-s, x-z) q(s, z, y) dz dy
\end{equation}
with function $q(s, z, y)$   be determined below.
We want
$$
\frac{\p }{\p t} p (t, x, y)= \LL^\kappa p(t, \cdot, y) (x) =\LL^{\kappa (x)} p  (t, \cdot, y) (x).
$$
Formally,
\begin{eqnarray*}
\frac{\p}{\partial t} p (t,x,y)
&=&  \LL^{\kappa(y)} p_y(t,x-y)+q(t,x,y)    +\int^t_0 \int_{\R^d}  \partial_t p_z(t-s,x-z)  q(s,z,y) d zd s \\
&=&  \LL^{\kappa(y)} p_y(t,x-y) + q(t,x,y)   +\int^t_0
\int_{\R^d}  \LL^{\kappa (z)}  p_z(t-s,x-z) q(s,z,y) d zd s  ,
\end{eqnarray*}
while
\begin{eqnarray*}
  \LL^{\kappa(x)} p (t,x,y)
&=&   \LL^{\kappa(x)} p_y(t,x-y) +\int^t_0
\int_{\R^d}  \LL^{\kappa(x)} p_z(t-s,x-z)  q(s,z,y)
dzd s \\
&=&    \LL^{\kappa(x)} p_y(t,x-y) +\int^t_0
\int_{\R^d}  q_0(t-s,x, z)  q(s,z,y)
dzd s   \\
&& +\int_{\R^d}  \LL^{\kappa(z)} p_z(t-s,x-z)  q(s,z,y)
dzd s,
\end{eqnarray*}
where
 $$ q_0(t, x, z)= (\LL^{\kappa(x)}-\LL^{\kappa(z)}) p_z (t,x-z).
$$
It follows from \eqref{e:3.12} that $q(t, x, y)$ should satisfy
\begin{equation} \label{e:3.14}
  q(t,x,y)=q_0(t,x,y)+\int^t_0 \int_{\R^d}q_0(t-s,x,z)  q(s,z,y)  dz d s.
\end{equation}
 Thus for the construction and  the upper bound heat kernel  estimates of $p(t, x, y)$,
the main task is to solve $q(t, x, y)$,  and to make the above argument rigorous.
We use Picard's iteration to solve \eqref{e:3.14}.
For $n\geq 1$, define
\begin{equation} \label{e:3.15}
  q_n(t,x,y)= \int^t_0 \int_{\R^d}q_0(t-s,x,z)  q_{n-1}(s,z,y)  dz d s.
\end{equation}
Then it can be shown that
\begin{equation} \label{e:3.16}
q(t, x, y):= \sum_{n=0}^\infty  q_n(t,x,y)
\end{equation}
 converges absolutely and locally uniformly
on $(0, 1]\times \R^d \times \R^d$. Moreover, $q(t, x, y)$ is jointly continuous in $(t, x, y)$ and has the following
upper bound estimate
$$
| q(t, x, y) | \leq C \left( \varrho^\beta_0 + \varrho^0_\beta \right) (t, x-y) ,
$$
where
$$ \varrho^\beta_\gamma (t, x)  := \frac{t^{\gamma /\alpha} (|x|^\beta \wedge 1)}{(t^{1/\alpha} + |x|)^{d+\alpha}}.
$$
We then need to address the following issues.
\begin{description}
\item{(i)} Show that $p(t, x, y)$ constructed through \eqref{e:3.13} and \eqref{e:3.16} is non-negative,
has the property $\int_{\R^d} p(t, x, y) dy=1$ and satisfies the Chapman-Kolmogorov equation.

\item{(ii)} The kernel $p(t, x, y)$ has the claimed two-sided estimates, and derivative estimates.

\item{(iii)} Uniqueness of $p(t, x, y)$.
\end{description}

This  requires  
detailed studies
on the kernel $p^\kappa_\alpha(t, x-y)$ for the symmetric L\'evy process
with L\'evy measure $\frac{\kappa (z)}{|z|^{d+\alpha}}dz $, including its fractional derivative estimates,
 and its continuous dependence on $\kappa (z)$, which will be outlined in the next two subsections.

\subsection{Upper bound estimates}

Key observation:  For any symmetric function $\kappa (z)$ with $\kappa_0\leq \kappa (z) \leq \kappa_1$,
let $\wh \kappa (z):= \kappa (z)-\frac{\kappa_0}2$.
Since the L\'evy process with L\'evy measure
$\frac{\kappa (z)}{|z|^{d+\alpha}}dz $ can be decomposed as the independent sum of L\'evy processes
having respectively L\'evy measures 
$\frac{\wh \kappa (z)}{|z|^{d+\alpha}}dz $ and $\frac{\kappa_0/2}{|z|^{d+\alpha}}dz$,
we have
$$
p^{\kappa(z)}_\alpha (t, x ) = \int_{\R^d} p^{\kappa_0/2}_\alpha (t, x-y) p^{\wh\kappa(z)}_\alpha (t, y) dy.
$$
 Thus the gradient and fractional derivative estimates on $p^{\kappa(z)}_\alpha (t, x )$ can be obtained
 from those on $p^{\kappa_0/2} (t, x ) $.
On the other hand, it follows from \cite{CK1} that there is a constant $c=c(d, \kappa_0, \kappa_1)\geq 1$ so that 
\begin{equation}\label{e:3.16b}
c^{-1} \varrho^0_{\alpha} (t, x) \leq p_\alpha^{\kappa (z)} (t, x) \leq  c\varrho^0_{\alpha} (t, x) 
\quad \hbox{for all } t>0 \hbox{ and } x\in \R^d.
\end{equation}

 First one can establish that for $\gamma_1, \gamma_2, \beta_1, \beta_2 \geq 0$
with $\beta_1 + \gamma_1>0$ and $\beta_2+\gamma_2>0$,
 \begin{eqnarray}\label{e:3.17}
&&  \int_0^t \int_{\R^d} \varrho^{\beta_1}_{\gamma_1} (t-s, x-z) \varrho^{\beta_2}_{\gamma_2} (s, z) dz ds
\nonumber \\
 &\leq & {\cal B}(\tfrac{\gamma_1+\beta_1}{\alpha } ,  \tfrac{\gamma_2+\beta_2}{\alpha } )
 \left( \varrho^0_{\gamma_1+\gamma_2+\beta_1+\beta_2}+ \varrho^{\beta_1}_{\gamma_1+\gamma_2+\beta_2}
 +\varrho^{\beta_2}_{\gamma_1+\gamma_2+\beta_1} \right) (t, x),
 \end{eqnarray}
 where ${\cal B}$ denotes   the usual $\beta$-function.

 Next we establish the continuous dependence of $p^{\kappa(z)}_\alpha (t, y)$ on the symmetric function $\kappa (z)$.
 Let $\kappa (z)$ and $\wt \kappa (z)$ be two symmetric functions that are bounded between $\kappa_0$ and $\kappa_1$.
 Then for every  $0<\gamma <\alpha/4$, there is a constant $c>0$ so that the following estimates hold for all $t\in (0, 1]$ and $x\in \R^d$,
\begin{eqnarray}
 |p^{\kappa(z)}_\alpha(t, x) - p^{\wt\kappa(z)}_\alpha(t, x) | &\leq& c \| \kappa -\wt \kappa \|_\infty \, \varrho^0_\alpha  (t, x), 
  \label{e:3.18} \\
 |\nabla_x p^{\kappa(z)}_\alpha(t, x) - \nabla_x p^{\wt\kappa(z)}_\alpha(t, x) |
 &\leq&  c \| \kappa -\wt \kappa \|_\infty  t^{-1/\alpha} \, \varrho^0_\alpha (t, x),
  \label{e:3.19} \\
 \int_{\R^d}    |\delta_{p^{\kappa}_\alpha } (t, x; z) - \delta_{p^{\wt \kappa}_\alpha } (t, x; z) |\frac{dz}{|z|^{d+\alpha}}  
 &\leq&   c \| \kappa -\wt \kappa \|_\infty \, \varrho^0_\alpha  (t, x). \label{e:3.20}
  \end{eqnarray}
Here $\| \kappa -\wt \kappa \|_\infty:= \sup_{x\in \R^d} | \kappa (z)-\wt \kappa (z)|$
and $\delta_f (t, x; z):=  f(t, x+z)+f (t, x-z) -2f(t, x)$.  
The above estimates are established in \cite[Theorem 2.5]{CZ}, but with an extra term on their right hand sides. For example, \eqref{e:3.18} corresponds to 
\cite[(2.30)]{CZ} where the estimate is
\begin{equation}\label{e:3.23a}
|p^{\kappa(z)}_\alpha(t, x) - p^{\wt\kappa(z)}_\alpha(t, x) |  \leq  c \| \kappa -\wt \kappa \|_\infty  \left( \varrho^0_\alpha + \varrho^\gamma_{\alpha-\gamma} \right) (t, x).
\end{equation}
 We take this opportunity to fill a gap in the proof of \cite[(3.20)]{CZ}.
The proof there
works only for $|x| \geq t^{1/\alpha}$ and $t\in (0, 1]$, as in this case,
by \cite[(2.2)]{CZ},
\begin{eqnarray*}
 \int_0^t \int_{\R^d} \varrho_0^\gamma (t-s, x-y) \varrho^0_{\alpha-\gamma} (s, x) dy ds
&\leq&  c_1 \varrho^0_{\alpha-\gamma} (t, x) \int_0^t \int_{\R^d} \varrho^\gamma_0(s,y)dy ds\\
&\leq & c_2  \varrho^0_{\alpha -\gamma } (t, x) t^{\gamma/\alpha} 
= c_2  \varrho^0_{\alpha} (t, x),
\end{eqnarray*}
which gives \eqref{e:3.23a} by line 8 on p.284 of \cite{CZ}.
On the other hand,  one deduces by the inverse Fourier transform that
$$
\sup_{y\in \R^d} \left| p_\alpha^{\kappa (z)} (t, y)-p_\alpha^{\wt \kappa (z)} (t, y)\right|
\leq (2\pi)^d \int_{\R^d} | e^{-t \psi_\kappa (\xi)} - e^{-t \psi_{\wt \kappa}  (\xi)} | d\xi 
\leq  c_3 \| \kappa -\wt \kappa \|_\infty  t^{-d/\alpha} .
$$ 
Thus when $|x|\leq t^{1/\alpha}$, 
$\left| p_\alpha^{\kappa (z)} (t, x)-p_\alpha^{\wt \kappa (z)} (t, x)\right| \leq c_3 \| \kappa -\wt \kappa \|_\infty  t^{-d/\alpha} 
\leq c_4 \| \kappa -\wt \kappa \|_\infty  \varrho^0_{\alpha} (t, x)$.
 
In fact, by a slight modification of the original proof given in \cite{CZ} for \eqref{e:3.23a}, we can 
get estimate \eqref{e:3.18}.
Indeed, by the symmetry of $\LL^{\kappa (z)}$ and $\LL^{\wt \kappa (z)}$,  
\begin{eqnarray*}
&& p^{\kappa (z)}_\alpha(t,x)-p^{\tilde\kappa (z)}_\alpha(t,x) \\
&=& \int^t_0 \frac{d}{ds} \left( \int_{\R^d}
p^{\kappa (z)}_\alpha(s,y) p^{\tilde\kappa (z)}_\alpha(t-s,x-y) dy\right) ds  \\
&=& \int^t_0 \left(\int_{\R^d} \left( \sL^{\kappa (z)}_\alpha p^{\kappa (z)}_\alpha(s,\cdot ) (y) 
 p^{\tilde\kappa (z)}_\alpha(t-s,x-y)
-p^{\kappa (z)}_\alpha(s,y) \sL^{\tilde\kappa (z)}_\alpha p^{\tilde\kappa (z)}_\alpha(t-s, \cdot) (x-y) \right)dy \right) ds\\
&=& \int^{t/2}_0 \left( \int_{\R^d}p^{\kappa (z)}_\alpha(s,y)
\left(\sL^{\kappa (z)}_\alpha-\sL^{\tilde\kappa (z)}_\alpha \right) p^{\tilde\kappa (z)}_\alpha(t-s, \cdot )
(x-y)dy \right) ds\\
&& +\int^t_{t/2} \left( \int_{\R^d}  p^{\tilde\kappa (z)}_\alpha(t-s,x-y) \left(\sL^{\kappa (z)}_\alpha-\sL^{\tilde\kappa (z)}_\alpha \right)
p^{\kappa (z)}_\alpha(s, \cdot)(y) dy \right) ds.
\end{eqnarray*}
Hence by  \eqref{e:3.16b} and  
\cite[(2.28)]{CZ}, 
\begin{eqnarray*}
|p^{\kappa (z)}_\alpha(t,x)-p^{\tilde\kappa (z)}_\alpha(t,x)|
&\leq& c\|\kappa -\tilde\kappa\|_\infty\int^{t/2}_0\!\int_{\R^d}\varrho^0_\alpha(s,y)\varrho^0_0(t-s,x-y)d yd s\\
&& +  c\|\kappa-\tilde\kappa\|_\infty\int^t_{t/2}\!\int_{\R^d}\varrho^0_0(s,y)\varrho^0_\alpha(t-s,x-y)d yd s\\
&\leq & \frac{c\|\kappa-\tilde\kappa\|_\infty }{t} \int^t_0 \!\int_{\R^d}\varrho^0_\alpha (s,y)\varrho^0_\alpha(t-s,x-y)d yd s \\
&\leq & c\|\kappa-\tilde\kappa\|_\infty  \, \varrho^0_\alpha (t,x).
\end{eqnarray*}
The same proof as that for \cite[Theorem 2.5]{CZ} but using \eqref{e:3.18} instead of
\eqref{e:3.23a} then gives  \eqref{e:3.19}-\eqref{e:3.20}.
 
Since
$$
\LL^{\kappa (z)} f(x)= {\rm p.v.} \int_{\R^d} \left( f(x+z)-f(x) \right) \frac{\kappa (z)}{|z|^{d +\alpha}} dz
=  \frac12 \int_{\R^d} \delta_f (x; z) \frac{\kappa (z)}{|z|^{d+\alpha}} dz, 
$$
estimate \eqref{e:3.20} implies that
 $$
 | \LL^{\kappa (z)} p^{\kappa(z)}_\alpha (t, x) - \LL^{\wt \kappa (z)} p^{\wt\kappa(z)}_\alpha (t, x)  |
 \leq  c \| \kappa -\wt \kappa \|_\infty  
\, \varrho^0_\alpha (t, x).
 $$
 From these estimates, one can establish the first part ((i)-(iv)) of the Theorem \ref{T:3.1}  as well as
 \begin{equation}\label{e:3.21}
 p(t, x, y) \geq c t^{-d/\alpha} \qquad \hbox{for } t\in (0, 1] \hbox{ and }
 |x-y| \leq 3 t^{1/\alpha}.
 \end{equation}

 \subsection{Lower bound estimates}

    The upper bound estimates in Theorem \ref{T:3.1} are established by using analytic method, while
 the lower bound estimate in Theorem \ref{T:3.1} are obtained mainly by probabilistic argument.

 From (i)-(iv) of Theorem \ref{T:3.1}, we see  that $P_t f (x):= \int_{\R^d} p(t, x, y) f(y)dy$ is a Feller semigroup.
Hence, it determines a Feller process $(\Omega,{\cal F}, (\P_x)_{x\in\R^d}, (X_t)_{t\geq 0})$.

We first claim the following.

\begin{thm}\label{T:3.2}
Let  ${\cal F}_t:=\sigma\{X_s, s\leq t\}$. Then for each $x\in \R^d$ and every $f\in C^2_b(\R^d)$, under $\P_x$,
\begin{equation}\label{ERY1}
M^f_t:=f(X_t)-f(X_0)-\int^t_0\LL f(X_s)d s  \ \hbox{   is an ${\cal F}_t$-martingale}.
\end{equation}
In other words, $\P_x$ solves the martingale problem for $(\LL,
C^2_b (\R^d))$. Thus $\P_x$ in particular solves the martingale problem
for $(\LL, C^\infty_c (\R^d))$.
\end{thm}

\noindent{\bf Sketch of Proof.}
For $f\in C^2_b(\R^d)$,  define $u(t, x)=f (x) + \int_0^t P_s \LL f(x) ds$.
Then we have by \eqref{eq15} in Theorem \ref{T:3.1} that
$$
\LL u(t, x) = \LL f (x)+ \int_0^t \LL P_s \LL f(x) ds = \LL f (x)+ \int_0^t \partial_s (P_s\LL f) (x)
= P_t \LL f (x) = \partial_t u(t, x).
$$
Since $P_t f$ also satisfies the equation $\partial_t P_t f =\LL (P_t f)$ with $P_0f =f$, we have
\begin{equation}\label{e:3.23}
P_t f (x) =f (x) + \int_0^t P_s \LL f(x) ds.
\end{equation}
The desired property \eqref{ERY1} now follows from \eqref{e:3.23} and the Markov property of $X$.
\qed

Theorem \ref{T:3.2} allows us to derive a L\'evy system of $X$ by following an approach from \cite{CKS1}.
It is easy to see from \eqref{ERY1} that $X_t=(X^{1}_t, \dots, X^{d}_t) $ is a semi-martingale.
For any $f\in C^\infty_c(\R^d)$, we have by It\^o's formula that
\begin{equation}\label{e:ito}
f(X_t)-f(X_0)=\sum^d_{i=1}\int^t_0{\partial}_if(X_{s-})d X^{i}_s +\sum_{s\leq t}\eta_s(f)  +\frac12 \gamma_t(f),
\end{equation}
where
\begin{equation}\label{e:ito2}
\eta_s(f)=f(X_s)-f(X_{s-})-\sum^d_{i=1}{\partial}_if(X_{s-})(X^{i}_s-X^{i}_{s-})
\end{equation}
and
\begin{equation}\label{e:ito3} 
\gamma_t(f)
=\sum^d_{i, j=1}\int^t_0{\partial}_i{\partial}_jf(X_{s-})d\langle X^{i, c}, X^{j,c}\rangle_s.
\end{equation}
Here $X^{i, c}$ is the continuous local martingale part of the semimartingale $X^i$.

Now suppose that $A$ and $B$ are two bounded closed subsets of $\R^d$
having a positive distance from each other. Let $f\in C^\infty_c(\R^d)$ with
$f=0$ on $A$ and $f=1$ on $B$.  Let $M^f$ be defined as in \eqref{ERY1}.
Clearly,
$N^f_t:=\int^t_0\1_A(X_{s-})d M^f_s$ is a martingale.
Define
\begin{equation}
 J(x, y)= k(x, y-x)/|y-x|^{d+\alpha},
\end{equation}\label{e:4.28}
so $\LL$ can be rewritten as
\begin{equation}\label{e:L2}
\LL f(x)= \lim_{\eps\to 0} \int_{\{|y-x|>\eps\}}
(f(y)-f(x)) J(x, y) d y.
\end{equation}
We get by  \eqref{ERY1}--\eqref{e:ito3} and  \eqref{e:L2},
\begin{align*}
N^f_t&=\sum_{s\leq t}\1_A(X_{s-})(f(X_s)-f(X_{s-})) -\int^t_0\1_A(X_s)\LL f(X_s)d s\\
&=\sum_{s\leq t}\1_A(X_{s-})f(X_s)-\int^t_0\1_A(X_s)\int_{\R^d}f(y)J(X_s,y)d yd s.
\end{align*}
By taking a sequence of functions $f_n\in C^\infty_c(\R^d)$ with
$f_n=0$ on $A$,  $f_n=1$ on $B$ and $f_n\downarrow \1_B$, we get
that, for any $x\in \R^d$,
$$
\sum_{s\leq t}\1_A(X_{s-})\1_B(X_s) -\int^t_0{\bf1}_A(X_s)\int_BJ(X_s, y)d yd s
$$
is a martingale with respect to $\P_x$. Thus,
$$
\E_x\left[ \sum_{s\leq t}\1_A(X_{s-})\1_B(X_s)\right]=
\E_x\left[\int^t_0\int_{\R^d} \1_A(X_s)\1_B(y)J(X_s, y)d yd s\right].
$$
Using this and a routine measure theoretic argument, we get
$$
\E_x\left[ \sum_{s\leq t}f(X_{s-}, X_s) \right]
=\E_x\left[\int^t_0\int_{\R^d}f(X_s, y)J(X_s, y)d yd s\right]
$$
for any non-negative measurable function $f$ on $\R^d\times \R^d$
vanishing on $\{(x, y)\in \R^d\times \R^d: x=y\}$. Finally,
following the same arguments as in \cite[Lemma 4.7]{CK1} and
\cite[Appendix A]{CK2}, we get

\begin{thm}\label{T:l2}
$X$ has  a L\'evy system $(J, t)$ with $J$ given by \eqref{e:4.28};  that is, for any
$x\in \R^d$ and any non-negative measurable function $f$ on $\R_+
\times \R^d\times \R^d$ vanishing on $\{(s, x, y)\in \R_+ \times
\R^d\times \R^d: x=y\}$ and $({\cal F}_t)$-stopping time $T$,
\begin{equation}\label{e:ls4xb}
\E_x \left[\sum_{s\leq T} f(s,X_{s-}, X_s) \right]= \E_x \left[
\int_0^T \left( \int_{\R^d} f(s,X_s, y) J(X_s, y)d y \right) d s
\right].
 \end{equation}
\end{thm}

For a set $K\subset\R^d$, denote
$$
\sigma_K:=\inf\{t\geq 0: X_t\in K \},\ \ \tau_K:=\inf\{t\geq 0: X_t\notin K\}.
$$
Let $B(x,r)$ be the ball with radius $r$ and center $x$. We need the following lemma (see \cite{BL02a, CK1}).
 
\begin{lemma}\label{Le56}
For each $\gamma\in(0,1)$, there exists $R_0>0$ such that for every $R>R_0$ and $r\in(0,1)$,
\begin{align}
\P_x(\tau_{B(x,Rr)} \leq r^\alpha)\leq \gamma.\label{ERY3}
\end{align}
\end{lemma}

\pf
Without loss of generality, we assume that $x=0$.
Given $f\in C^2_b(\R^d)$ with $f(0)=0$ and $f(x)=1$ for $|x|\geq 1$,  we set
$$
f_r(x):=f(x/r),\ \ r>0.
$$
By the definition of $f_r$, we have
\begin{align}
 \P_0(\tau_{B(0,Rr)}\leq r^\alpha)
\leq \E_0 \left[ f_{Rr} ( X_{\tau_{B(0,Rr)}\wedge r^{\alpha}}) \right]
\stackrel{(\ref{ERY1})}{=}\E_0\left(\int^{\tau_{B(0,Rr)}\wedge r^{\alpha}}_0\LL f_{Rr}(X_s)d s\right).\label{ERY2}
\end{align}
On the other hand, by the definition of $\LL $, we have for $\lambda>0$,
\begin{align*}
|\LL f_{Rr}(x)|&=\frac{1}{2}\left|\int_{\R^d}(f_{Rr}(x+z)+f_{Rr}(x-z)-2f_{Rr}(x))\kappa(x,z)|z|^{-d-\alpha}d z\right|\\
&\leq\frac{\kappa_1\|\nabla^2 f_{Rr}\|_\infty}{2}\int_{|z|\leq\lambda r}|z|^{2-d-\alpha}d z+2\kappa_1\|f_{Rr}\|_\infty\int_{|z|\geq\lambda r}|z|^{-d-\alpha}d z\\
&=\kappa_1\frac{\|\nabla^2 f\|_\infty}{(Rr)^2}\frac{(\lambda r)^{2-\alpha}}{2(2-\alpha)}s_1+2\kappa_1\|f\|_\infty\frac{(\lambda r)^{-\alpha}}{\alpha}s_1\\
&=\kappa_1s_1\left(\frac{\|\nabla^2 f\|_\infty}{R^2}\frac{\lambda^{2-\alpha}}{2(2-\alpha)}+2\|f\|_\infty
\frac{\lambda^{-\alpha}}{\alpha}\right) r^{-\alpha},
\end{align*}
where $s_1$ is the sphere area of the unit ball. Substituting this into (\ref{ERY2}),
we get
$$
\P_0(\tau_{B(0,Rr)}\leq r^\alpha) \leq
\kappa_1s_1\left(\frac{\|\nabla^2 f\|_\infty}{R^2}\frac{\lambda^{2(2-\alpha)}}{2-\alpha}+2\|f\|_\infty
\frac{\lambda^{-\alpha}}{\alpha}\right) .
$$
Choosing
first $\lambda$ large enough and then $R$ large enough
yield the desired estimate.
\qed

\medskip

We can now proceed to establish the lower bound heat kernel estimate \eqref{e:3.10}.
By Lemma \ref{Le56},  there is a constant $\lambda\in (0,\tfrac{1}{2})$ such that for all $t\in(0,1)$,
\begin{equation}\label{e:3.32}
\P_x(\tau_{B(x,t^{1/\alpha}/2)}> \lambda t)\geq\tfrac{1}{2}.
\end{equation}
In view of the estimate \eqref{e:3.21},  it remains to consider the case that $|x-y|>3t^{1/\alpha}$.
Using \eqref{e:3.32} and the L\'evy system of $X$,
 \begin{eqnarray*}
&& \P_x(X_{\lambda t}\in B(y,   t^{1/\alpha})) \nonumber \\
&\geq& \P_x\left( X \hbox{ hits } B(y,  t^{1/\alpha}/2) \hbox{ before }
\lambda  t \hbox{ and then travels less than } \right. \nonumber \\
&& \hskip 1.0truein \left. \hbox{distance }
t^{1/\alpha}/2 \hbox{ for at least }
 \lambda t \hbox{ units of time} \right) \nonumber \\
&\geq & \P_x(\sigma_{B(y, t^{1/\alpha}/2)}<\lambda t)\inf_{z\in B(y, t^{1/\alpha}/2)}\P_z(\tau_{B(z,
t^{1/\alpha}/2)}> \lambda t) \nonumber \\
&\geq & c_1 \P_x (X_{(\lambda t) \wedge \tau_{B(x,t^{1/\alpha})}}\in B(y,t^{1/\alpha}/2) ) \\
&=& \E_x\int_0^{ (\lambda t) \wedge\tau_{B(x,t^{1/\alpha})}}\int_{B(y,t^{1/\alpha}/2)}
J(X_s, u)\,du\,ds \\
 &\geq & c_2\E_x \left[ (\lambda t) \wedge
\tau_{B(x,t^{1/\alpha})}\right]\int_{B(y,t^{1/\alpha}/2)}
 \frac{1}{|x-y|^{d+\alpha}}   du \\
&\geq&  c_3\, \frac{t^{(d+\alpha)/\alpha}}{|x-y|^{d+\alpha}}   .
\end{eqnarray*}
 Thus
\begin{eqnarray*}
p(t,x,y)&\geq& \int_{B(y,t^{1/\alpha})} p(\lambda t,x,z) p((1-\lambda) t, z, y)\,dz\\
&\geq&  \P_x(X_{\lambda t}\in B(y,t^{1/\alpha})) \, \inf_{z\in B(y,t^{1/\alpha})}  p((1-\lambda t,z,y) \\
&\geq&  c_4 t^{-d/\alpha} \, t^{(d+\alpha)/\alpha}
  \frac{1}{|x-y|^{d+\alpha}}  \\
&=& \frac{c_4 t}{|x-y|^{d+\alpha}} .
\end{eqnarray*}
This proves that
$$
p(t, x, y) \geq c \left(t^{-d/\alpha}\wedge \frac{t}{|x-y|^{d+\alpha}}\right)
\quad \hbox{for every } x, y \in \R^d \hbox{ and } t\leq 1.
$$

\subsection{Strong stability}

In real world applications and modeling, state-dependent parameter  $\kappa (x, z)$ of \eqref{e:3.1} is an
approximation of real data. So  a natural question
is how reliable the conclusion is when using such an approximation.
The following strong stability result is recently obtained in \cite{CZ2}.

\begin{thm}\label{T:3.5}
Suppose $\beta \in (0, \alpha/4]$, and $\kappa$ and $\wt  \kappa$ are two functions satisfying  \eqref{e:2.6} and  \eqref{e:2.7}.
Denote the corresponding fundamental solution by $p^\kappa  (t,x,y) $ and $p^{\wt  \kappa}  (t,x,y)$, respectively.
 Then for every
$\gamma \in (0, \beta)$ and $\eta\in (0, 1)$, there exists a constant
$C=C(d, \alpha, \beta, \kappa_0, \kappa_1, \kappa_2, \gamma, \eta)>0$
so that for all $t\in (0, 1]$ and $x, y\in \R^d$,
\begin{equation}\label{e:3.33}
| p^\kappa_{\alpha} (t,x,y) -p^{\wt  \kappa}_{\alpha} (t,x,y) |
\leq C  \| \kappa -  \wt  \kappa \|^{1-\eta}_\infty \,   \left( 1+t^{-\gamma/\alpha} (|x-y|^\gamma \wedge 1)\right)
\, \frac{t} {(t^{1/\alpha} + |x-y|)^{d+\alpha}}.
\end{equation}
Here $\|  \kappa -  \wt  \kappa \|_\infty :=\sup_{x, z\in \R^d} | \kappa (x, z)-\wt  \kappa (x, z)|$.
 \end{thm}

 Observe that by \eqref{eq16} and \eqref{e:3.10},  
 the term $\frac{t} {(t^{1/\alpha} + |x-y|)^{d+\alpha}}$
 in \eqref{e:3.33} is comparable to
$p^\kappa_\alpha(t, x, y)$ and to $p^{\wt  \kappa}(t, x, y)$.
So the error bound \eqref{e:3.33} is also a relative error bound, which is good even in the region when $|x-y|$ is large.

Let $\{P_t^\kappa; t\geq 0\}$ and $\{P_t^{\wt  \kappa}; t\geq 0\}$
be the semigroups generated by
$\sL^{\kappa} $ and $\sL^{\wt  \kappa} $, respectively.
For $p\geq 1$, denote by
$\| P^\kappa_t - P^{\wt  \kappa}_t\|_{p, p} $ the operator norm
of $P^\kappa_t - P^{\wt  \kappa}_t$ in Banach space
$L^p (\R^d; dx)$. 

\begin{cor}\label{C:3.6}
Suppose $\beta \in (0, \alpha/4]$, and $\kappa$ and $\wt  \kappa$ are two functions satisfying  \eqref{e:2.6} and  \eqref{e:2.7}.
 Then for every
$\gamma \in (0, \beta)$ and $\eta\in (0, 1)$, there exists a constant
 $C=C (d, \alpha, \beta, \kappa_0, \kappa_1, \kappa_2, \gamma, \eta)>0$
so that for every $p\geq 1$ and $t\in (0, 1]$,
\begin{equation}\label{e:3.34}
\| P^\kappa_t - P^{\wt  \kappa}_t\|_{p, p} \leq
 C  t^{-\gamma /\alpha} \| \kappa -  \wt  \kappa \|^{1-\eta}_\infty .
\end{equation}
\end{cor}

  Theorem \ref{T:3.5} is derived by estimating each $|q^\kappa_n(t, x, y)-q_n^{\wt \kappa}(t, x, y)|$
    for  $q^\kappa_n(t, x, y)$ and $q_n^{\wt \kappa}(t, x, y)$ of \eqref{e:3.15}. Corollary \ref{C:3.6}
    is a direct consequence of Theorem \ref{T:3.5}. 

For uniformly elliptic divergence form operators $\sL$ and $\wt  \sL$ on $\R^d$,
pointwise estimate on $|p(t, x, y)-\wt  p(t, x, y)|$
and the $L^p$-operator norm estimates on
 $P_t -\wt  P_t$
are obtained in Chen, Hu, Qian and Zheng \cite{CHQZ}
in terms of the local $L^2$-distance between the diffusion
matrix of $\sL$ and   $\wt  \sL$.
Recently, Bass and Ren \cite{BR} obtained strong stability result
for symmetric $\alpha$-stable-like non-local operators
of \eqref{e:3.3}, with error bound expressed
in terms of the $L^q$-norm
on the function $c(x):=\sup_{y\in \R^d}|c(x, y)-\wt  c(x, y)|$.

\subsection{Applications to SDE driven by stable processes}
 
 Suppose that $\sigma(x)=(\sigma_{ij}(x))_{1\leq i, j\leq d}$
is a bounded continuous $d\times d$-matrix-valued function
on $\R^d$ that is  non-degenerate
at every $x\in \R^d$,
and $Y$ is a (rotationally) symmetric $\alpha$-stable process
on $\R^d$ for some $0<\alpha <2$. It is shown in Bass and Chen
\cite[Theorem 7.1]{BC} that for every $x\in \R^d$, SDE
\begin{equation}\label{e:1.18}
d X_t = \sigma(X_{t-}) d Y_t, \qquad X_0=x,
\end{equation}
has a unique weak solution.
(Although in \cite{BC} it is assumed $d\geq 2$, the
 results there are valid for $d=1$ as well.)
The family of these weak solutions
forms a strong Markov process
$\{X, \P_x, x\in \R^d\}$.
Using It\^o's
formula, one deduces (see the display above (7.2) in
\cite{BC}) that $X$ has   generator
\begin{equation}\label{e:1.19}
\LL f(x) = {\rm p.v.} \int_{\R^d} \left( f(x+\sigma(x)y)-f(x)\right)
 \frac{{\cal A}(d, -\alpha) } {|y|^{d+\alpha}} d y.
 \end{equation}
 A change of variable formula $z=\sigma(x) y$ yields
\begin{equation}\label{e:1.20}
\LL f(x) = {\rm p.v.} \int_{\R^d} \left( f(x+z)-f(x)\right)
\frac{\kappa (x, z)} {|z|^{d+\alpha}} d z ,
\end{equation}
where
\begin{equation}\label{e:1.21}
 \kappa (x, z)= \frac{{\cal A}(d, -\alpha)}{|{\rm det} \sigma(x)|} 
\left( \frac{|z|}{|\sigma(x)^{-1}z|}\right)^{d+\alpha}.
\end{equation}
Here ${\rm det}(\sigma(x))$ is the determinant of the matrix $\sigma(x)$
and $\sigma(x)^{-1}$ is the inverse of $\sigma(x)$.
As an application of Theorem \ref{T:3.1}, we have

\begin{cor}  {\rm (\cite[Corollary 1.2]{BC})} \label{C:1.4}
Suppose that $\sigma(x)=(\sigma_{ij}(x))$ is  a $d\times d$ matrix-valued function on $\R^d$ 
such that there are positive constants
$\lambda_0$,  $\lambda_1$, $\lambda_2$ and $\beta\in (0, 1)$   so that
\begin{equation}\label{e:3.40}
\lambda_0 {\bf I}_{d\times d} \leq \sigma(x)\leq \lambda_1 {\bf I}_{d\times d}
\quad \hbox{for every } x\in \R^d,
\end{equation}
and  
\begin{equation}\label{e:3.41}
 |\sigma_{ij}(x)-\sigma_{ij}(y)| \leq \lambda_2 |x-y|^\beta
 \quad \hbox{for } 1\leq i, j\leq d.
\end{equation}
Then the strong Markov process $X$ formed by the unique weak solution to SDE \eqref{e:1.18} has a jointly continuous transition density function
$p(t, x, y)$ with respect to the Lebesgue measure on $\R^d$, and
there is a constant  $C\geq 1$
that depends only on $(d, \alpha, \beta, \lambda_0, \lambda_1)$ so that
$$
C^{-1}\,  \frac{t}{(t^{1/\alpha} + |x-y|)^{d+\alpha}}
\leq p(t, x, y)   \leq
C\, \frac{t}{(t^{1/\alpha} + |x-y|)^{d+\alpha}}
$$
for every $t\in (0, 1]$ and $x, y\in \R^d$.
Moreover, $p(t, x, y)$ enjoys all the properties
 stated in the conclusions
of Theorem \ref{T:3.1}
 with $\kappa_0= {\cal A}(d, -\alpha)\lambda_0^{d+\alpha} \lambda_1^{-d}$,
 $\kappa_1= {\cal A}(d, -\alpha) \lambda_0^{-d} \lambda_1^{d+\alpha}$
and $\kappa_2= \kappa_2 (d, \lambda_0, \lambda_1, \lambda_2)$.
\end{cor}

The following strong stability result for SDE \eqref{e:1.18} is a direct consequence
of Corollary \ref{C:3.6}
and \eqref{e:1.21}. 

\begin{cor}  \label{C:3.8}
Suppose that $\sigma(x)=(\sigma_{ij}(x))$ and $\wt \sigma(x)=(\wt \sigma_{ij}(x))$ are  $d\times d$ matrix-valued functions on $\R^d$
satisfying conditions \eqref{e:3.40} and \eqref{e:3.41}. Let $p(t, x, y)$ and $\wt p(t, x, y)$ 
be the transition density functions of the corresponding  strong Markov processes $X$ and $\wt X$ 
that solve SDE \eqref{e:1.18}, respectively. 
Then for every
$\gamma \in (0, \beta)$ and $\eta\in (0, 1)$, there exists a constant
$C=C(d, \alpha, \beta, \lambda_0, \lambda_1, \lambda_2, \gamma, \eta)>0$
so that for all $t\in (0, 1]$ and $x, y\in \R^d$,
\begin{equation}\label{e:3.42}
| p (t,x,y) -\wt p  (t,x,y) |
\leq C  \| \sigma-  \wt  \sigma\|^{1-\eta}_\infty \,   \left( 1+t^{-\gamma/\alpha} (|x-y|^\gamma \wedge 1)\right)
\, \frac{t} {(t^{1/\alpha} + |x-y|)^{d+\alpha}},
\end{equation}
where  $\| \sigma-  \wt  \sigma\|^{1-\eta}_\infty:=\sum_{i,j=1}^d \sup_{x, y\in \R^d} |\sigma_{ij}(x)-\wt \sigma_{ij}(x)|$.
\end{cor}

 \section{Diffusion with jumps}\label{S:4}

 In this section, we consider non-local operators that have both elliptic differential operator part
 and pure non-local part:
 \begin{equation}\label{4.1}
  \LL f(x):=\LL^{a} f(x)+b   \cdot \nabla f(x)+\LL^{\kappa} f(x),
 \end{equation}
 where
 \begin{eqnarray*}
  && \LL^{a} f(x):=\frac{1}{2}\sum_{i,j=1}^d a_{ij} (x) \p ^2_{ij} f(x),\qquad  b \cdot \nabla f(x):=\sum_{i=1}^{d}b_i( x)\p_if(x),\\
  && \LL^{\kappa} f(x):=\int_{\R^d}\left(f(x+z)-f(x)- \1_{\{|z|\leq 1\}}z\cdot\nabla f(x)\right)\frac{\kappa (x,z)}{|z|^{d+\alpha}}d z.
\end{eqnarray*}
Here $a (x):=(a_{ij} (x))_{1\leq i,j\leq d}$ is a $d\times d$-symmetric matrix-valued measurable function on $ \R^d$, $b(x):  \R^d\to\R^d$ and $\kappa( x,z):  \R^d\times\R^d\to\R$ are measurable functions, and $\alpha\in(0,2)$.

For convenience, we assume $d\geq 2$.
Throughout this section, we impose the  following assumptions on $a$ and $\kappa$:
\begin{description}
  \item{\rm {(\bf H$^a$)}} There are $c_1>0$ and $\beta\in(0,1)$ such that for any  $x,y\in\R^d$,
  \begin{equation}
    |a(x)-a(y)|\leq c_1|x-y|^{\beta},\label{eqa2}
  \end{equation}
  and for some $c_2\ge1$,
  \begin{equation}
    c_2^{-1} {\bf I}_{d\times d} \leq a(x)\leq c_2  {\bf I}_{d\times d}.\label{eqa1}
  \end{equation}

  \item{\rm {(\bf H$^\kappa$)}}  $\kappa(x,z)$ is a bounded measurable function and if $\alpha=1$, we require
  \begin{equation}\label{Sym} 
      \int_{r<|z|\leq R}\kappa( x,z)|z|^{-d-1}d z=0 \quad \hbox{for any }  0< r<R<\infty.
  \end{equation}
\end{description}

Note that when $\kappa (x, z)$ is a positive constant function,
$$ \LL = \LL^a +b \cdot \nabla + c \Delta^{\alpha/2}$$
for some constant $c>0$.
A function $f$ defined on $\R^d$ is said to be in Kato class $\mK_2$ if
$f\in L^1_{loc}(\R^d)$ and
\begin{align}\label{ET3}
 \lim_{\delta \to 0} \sup_{x\in \R^d} \int_0^\delta \int_{\R^d}  \frac{t^{1/2}|f(y)|}{(t^{1/2}+|x-y|)^{d+2}} dy dt =0.
\end{align}

Let $q(t, x , y)$ be the fundamental solution of $\{\sL^a; t\geq 0\}$; see Theorem 
\ref{T:4.3} below for more information.
Since $\sL$ can be viewed as a perturbation of $\sL^a $ by $\sL^{b, \kappa} :=b\cdot \nabla + \sL^\kappa$, heuristically the fundamental solution
(or heat kernel)
$p (t, x, y)$
of $\sL $ should satisfy the following
Duhamel's formula: for all $t>0$ and $x, y\in \R^d$,
\begin{equation}\label{eqdu}
p(t, x, y)=q(t, x, y)+\int_0^t\!\!\! \int_{\R^d} p(r, x,  z) \sL^{b, \kappa} 
q(t-r, \cdot ,  y) (z) dz dr
\end{equation}
or
\begin{equation}\label{eqdu0}
p(t,  x, y)=q(t, x, y)+\int_0^t\!\!\! \int_{\R^d} q(r, x, z) \sL^{b, \kappa}
p(t-r, \cdot , y) (z) d z dr.
\end{equation}

 The following is a special case of the main results in \cite{CHXZ}, where the corresponding results are also obtained for time-inhomogeneous operators.

 \begin{thm}{\rm (\cite[Theorem 1.1]{CHXZ})} \label{T:4.1}
Let $\alpha\in(0,2)$. Under {\bf (H$^a$)}, {\bf (H$^\kappa$)} and $b\in\mK_2$, 
there is a unique  continuous function $p(t,x;   y)$  that   satisfies \eqref{eqdu},
and
\begin{description}
  \item{\rm (i)} (Upper-bound estimate) For any $T>0$, there exist constants $C_0,\lambda_0>0$ such that  for $t\in (0, T]$ and $x, y\in \R^d$,
      \begin{equation}     \label{eqlpe}
        |p(t, x, y)|\leq C_0 \left( t^{-d/2} e^{-\lambda_0 |x-y|^2/t} + \frac{\|\kappa\|_\infty \, t}{(t^{1/2}+|x-y|)^{d+\alpha}} \right) .
       \end{equation}

  \item{\rm (ii)}  (C-K equation) For all $s, t>0$ and $x, y\in \R^d$, we have
      \begin{equation}\label{CK01}
        \int_{\R^d} p(s, x, y) p(t,  y, z) dy = p(s+t, x, y).
        \end{equation}
        
  \item{\rm (iii)}  (Gradient estimate) For any $T>0$, there exist constants $C_1,\lambda_1>0$ such that
  for $t\in (0, T]$ and $x, y\in \R^d$,
      \begin{equation} \label{eqpge}
        |\nabla_x p(t, x, y)|\leq C_1 t^{-1/2} \left( t^{-d/2} e^{-\lambda_1 |x-y|^2/t}
        + \frac{\|\kappa\|_\infty \, t}{(t^{1/2}+|x-y|)^{d+\alpha}} \right) .
      \end{equation}

  \item{\rm (iv)}  (Conservativeness) For any $  t>0$ and $x \in \mathbb{R}^d$,
    $  \int_{\R^d}p(t , x, y)  d  y=1$.

  \item{\rm (v)} (Generator) Define $P_t f (x)=\int_{\R^d} p(t, x, y) f(y) dy$. Then for any $f\in C_b^2(\R^d)$, we have
      \begin{equation}\label{eqge}
        P_{t}f(x)-f(x)=\int_0^t  P_s \LL f(x) ds .
      \end{equation}

  \item{\rm (vi)} (Continuity) For any bounded and uniformly continuous function $f$, $\lim_{t\to 0}\|P_{t}f-f\|_\infty=0$.
  \end{description}
\end{thm}

 Define $m_\kappa = \inf_{x\in \R^d} {\rm essinf}_{z\in \mathbb{R}^d} k( x,z)$.

\begin{thm}{\rm (\cite[Theorem 1.3]{CHXZ})} \label{T:4.2}
  If $\kappa$ is a bounded function satisfying $({\bf H}^{\kappa})$ and that for  each $x \in \mathbb{R}^d$,
  \begin{equation}\label{e:4.9}
    \kappa( x,z) \ge 0 \quad \hbox{for a.e.  }   z\in \mathbb{R}^d,
  \end{equation}
  then $p(t,x ,y) \ge 0$.  
    Furthermore, if $m_\kappa>0$, then
  for any $T > 0$, there are constants $C_1, \lambda_2 > 0$ such that
  for any $t\in (0, T]$ and $x, y\in \R^d$,
  \begin{equation}\label{Low}
    p(t,x ,y)\ge C_1 \left( t^{-d/2} e^{-\lambda_2 |x-y|^2/t} + 
    \frac{ m_\kappa \, t}
    {(t^{1/2}+|x-y|)^{d+\alpha}} \right) .
      \end{equation}
\end{thm}

 We have by Theorems \ref{T:4.1} and \ref{T:4.2} that when $\kappa \geq 0$, then there is a conservative Feller process
 $X= \{X_t, t\geq 0; \P_x, x\in \R^d\}$ having $p(t, x, y)$ as its transition density function with respect to the Lebesgue measure.
 It follows from \eqref{eqge} that $X$ is a solution to the martingale problem for $(\LL, C^2_b (\R^d))$.

When $a$ is the identity matrix, $b=0$ and $\kappa (x, z)$ is a positive constant,
$\LL=\Delta + c \Delta^{\alpha/2}$ for some positive constant $c>0$. In this case, the corresponding
Markov process $X$ is a symmetric L\'evy process that is the sum of a Brownian motion $W$
and an independent rotationally symmetric $\alpha$-stable process $Y$. Thus the heat kernel $p(t, x, y)$ for $\LL$ is
the convolution of the transition density function of $W$ and $Y$. In this case, its two-sided bounds can be obtained through
a direct calculation. Indeed such a computation is carried out in Song and Vondra\v cek \cite{SV}.

Symmetric diffusions with jumps corresponding to symmetric non-local operators on $\R^d$ with variable coefficients of the the following form  have been studied in \cite{CK3}:
\begin{equation}\label{e:4.11}
 \LL f (x) = \sum_{i,j=1}^d \frac{\partial}{\partial x_i} \left( a_{ij}(x) \frac{\partial f (x)}{\partial x_i} \right)
+ \lim_{\eps \to 0} \int_{|x-y|>\eps} (f(y)-f(x)) \frac{c(x, y)}{|x-y|^{d+\alpha}} dy,
\end{equation}
where $a (x):=(a_{ij} (x))_{1\leq i,j\leq d}$ is a $d\times d$-symmetric matrix-valued measurable function on $ \R^d$, 
 $c(x, y)$ is a symmetric measurable function on $ \R^d\times\R^d$  
that is bounded between two positive constants,  and $\alpha\in(0,2)$.
Clearly, when $a(x)$ is the identity matrix and $c(x, y)$ is a positive constant, the above non-local operator is $\Delta + c_0 \Delta^{\alpha/2}$ for some $c_0>0$. 
Among other things, 
it is established in Chen and Kumagai \cite{CK3} that the symmetric non-local operator
$\LL$ of \eqref{e:4.11} has  
 a jointly H\"older continuous heat kernel $p(t, x, y)$ and there are positive constants $c_i$, $1\leq i\leq 4$ so that
\begin{eqnarray} \label{e:4.12}
&& c_1 \left( t^{-d/2}\wedge t^{-d/\alpha}\right ) \wedge \left(t^{-d/2} e^{-c_2|x-y|^2/t} + t^{-d/\alpha} \wedge \frac{t}{|x-y|^{d+\alpha}} \right) \nonumber  \\
&\leq & p(t, x, y) \leq c_3  \left( t^{-d/2}\wedge t^{-d/\alpha}\right ) \wedge
\left(t^{-d/4} e^{-c_4|x-y|^2/t} + t^{-d/\alpha} \wedge \frac{t}{|x-y|^{d+\alpha}} \right)
\end{eqnarray}
for all $t>0$ and $x, y\in \R^d$. It is easy to see that for each fixed $T>0$, the two-sided estimates \eqref{e:4.12} on $(0, T]\times \R^d \times \R^d$
is equivalent to
 $$
  \wt c_1 \left( t^{-d/2} e^{-c_2 |x-y|^2/t} + \frac{  t}{(t^{1/2}+|x-y|)^{d+\alpha}} \right)
  \leq p(t, x, y) \leq \wt c_3 \left( t^{-d/2} e^{-c_4  |x-y|^2/t} + \frac{  t}{(t^{1/2}+|x-y|)^{d+\alpha}} \right) .
$$

\medskip

 When $a$ is the identity matrix and $b=0$, the results in Theorems \ref{T:4.1} and \ref{T:4.2} have been obtained recently in \cite{Wang}
for $\kappa (x, z)$ that is symmetric in $z$.

 \subsection{Approach}

 The approach in \cite{CHXZ} is to treat $\LL$ as $\LL^a$ under lower order perturbation $b\cdot \nabla + \LL^\kappa$,
 and thus one can construct the fundamental solution for $\LL$ from that of $\LL^a$ through Duhamel's formula.

 The following result is essentially known in literature; see \cite{Fr} (see also \cite[Theorem 2.3]{CHXZ}).

 \begin{thm}\label{T:4.3}
  Under {\bf (H$^a$)}, there exists a nonnegative continuous function $q(t, x, y) $, called the fundamental solution or heat kernel of $\LL^a$,
   with the following properties:
  \begin{description}
    \item{\rm (i)}  (Two-sided estimates) For any $T>0$, there exist constants $C,\lambda>0$ such that
    for $t\in (0, T]$ and $x, y\in \R^d$,
        \begin{align}
          C^{-1} t^{-d/2} e^{-\lambda^{-1}|x-y|^2/t} \leq q(t,x , y)\leq C t^{-d/2} e^{-\lambda |x-y|^2/t} .\label{ET41}
        \end{align}

    \item{\rm (ii)}  (Gradient estimate) For $j=1,2$ and $T>0$, there exist constants $C,\lambda>0$ such that
    for $t\in (0, T]$ and $x, y\in \R^d$,
      \begin{equation}
          |\nabla^j_x  q(t, x ,y)|\leq C t^{-(d+j)/2} e^{-\lambda |x-y|^2/t} .   \label{eq21}
        \end{equation}

    \item{\rm (iii)}  (H\"older estimate in $y$) For $j=0,1$, $\eta\in(0,\beta)$ and $T>0$, there exist constants $C,\lambda>0$ such that
    for $t\in (0, T]$,  $x, y, z\in \R^d$,
  \begin{equation}\label{eq202}
          |\nabla^j_x  q(t,x, y)-\nabla^j_x  q(t, x, z)|\leq C|y-z|^{\eta} \,  t^{-(d+j+\eta)/2} \left( e^{-\lambda |x-y|^2/t} +
          e^{-\lambda |x-z|^2/t} \right).
  \end{equation}
 Moreover, for bounded measurable $f:\R^d\to\R$, let $Q_t f(x):=\int_{\R^d} q(t,x, y)f(y)d y$. We have

    \item{\rm (iv)} (Continuity) For any bounded and uniformly continuous function $f$,
      $  \lim_{t\to 0}\|Q_t f-f\|_\infty=0$.

     \item{\rm (v)}  (C-K equation) For all $0\leq t<r<s<\infty$,  $Q_t Q_s=Q_{t+s}$.

    \item{\rm (vi)} (Conservativeness) For all $0\leq t<s<\infty$,  $Q_t1=1$.

   \item{\rm (vii)}  (Generator) For any $f\in C^2_b(\R^d)$, we have
        $$          Q_tf(x)-f(x)=\int_0^t Q_s \LL^{a} f(x) dr=\int_0^t \LL^{a} Q_s f(x) ds.
        $$
 \end{description}
\end{thm}

As mentioned earlier, it is
expected that the fundamental solution $p(t, x, y)$ of $\LL$ should satisfy Duhamel's formula \eqref{eqdu}. 
We construct $p(t, x, y)$ recursively. Let $p_0(t, x, y)=q(t, x, y)$, and define for $n\geq 1$,
$$
p_n(t,x, y):=\int_0^t \int_{\R^d}p_{n-1}(t-s ,x, z)\LL^{b,\kappa} q(s, \cdot , y) (z) dz ds  .
$$
Using Theorem \ref{T:4.3}, one can show that $p_n(t, x, y)$ is well defined and that
$\sum_{n=0}^\infty p_n(t, x, y)$ converges locally uniformly to some function $p(t, x, y)$,
and that $p(t, x, y)$ is the unique solution stated in Theorem \ref{T:4.1}.
The positivity \eqref{e:4.9} of Theorem \ref{T:4.2} can be established by using
Hille-Yosida-Ray theorem  and Courr\'ege's first theorem.

The Gaussian part in the lower bound estimate on $p(t, x, y)$ in Theorem \ref{T:4.2} is obtained from 
the near diagonal lower bound estimate on $p(t, x, y)$ and a chaining ball argument, while the pure jump part in the lower bound estimate on $p(t, x, y)$ is obtained 
by using a probabilistic argument through the L\'evy system, 
similar to that in Section 3.

 \subsection{Application to SDE}
  Let $\sigma (x)$ be  a $d\times d$-matrix valued function on $\R^d$ that is uniformly elliptic
 and bounded, and each entry $\sigma_{ij}$ is $\beta$-H\"older continuous on $\R^d$,
 $b\in \mK_2$ and $\wt \sigma$ a  bounded $d\times d$-matrix valued measurable function on $\R^d$.
 Suppose $X$ solves the following stochastic differential equation:
 $$
 dX_t = \sigma (X_t) dB_t  + b(X_t) dt + \wt \sigma(X_{t-}) dY_t,
 $$
 where $W$ is a Brownian motion on $\R^d$ and $Y$ is a rotationally symmetric $\alpha$-stable process
 on $\R^d$. By It\^o's formula,  the infinitesimal generator $\LL$ of $X$ is of the form $\LL^a+b\cdot \nabla + \LL^\kappa$
 with $a(x)= \sigma (x) \sigma (x)^{*}$ and
 $$
 \kappa (x, z)= \frac{{\cal A}(d, -\alpha)}{|{\rm det} \wt \sigma(x)|} \left( \frac{|z|}{|\wt \sigma(x)^{-1}z|}\right)^{d+\alpha}.
   $$
 So by Theorems \ref{T:4.1} and \ref{T:4.2}, $X$ has a transition density function $p(t, x, y)$ satisfying the properties
 there. If in addition, $\wt \sigma$ is uniformly elliptic, 
   then for any $T>0$,
 $$
       c_1 \left( t^{-d/2} e^{- \lambda_1 |x-y|^2/ t } + \frac{  t}{(t^{1/2}+|x-y|)^{d+\alpha}} \right)
      \leq p(t, x, y) \leq c_2    \left( t^{-d/2} e^{- \lambda_2 |x-y|^2 /t } + \frac{  t}{(t^{1/2}+|x-y|)^{d+\alpha}} \right)
$$
   for $t\in(0,T]$ 
 and $x, y\in \R^d$.

 \section{Other related work}\label{S:5}

 In this section, we briefly mention some other recent work on heat kernels of non-symmetric non-local operators.

 Using a perturbation argument, Bogdan and Jakubowski \cite{BJ}
 constructed {\it a particular} heat kernel (also called fundamental solution) $q^b(t, x, y)$
for operator $\sL^b:=\Delta^{\alpha/2}+b\cdot \nabla$ on $\R^d$, where $d\geq 1$, $\alpha \in (1, 2)$
and $b$ is a function on $\R^d$ that is in a suitable Kato class.
It is based on the following heuristics: $q^b(t, x, y)$ of $\sL^b$ can be related to
 the fundamental solution $p(t, x, y)$ of $\sL^0=\Delta^{\alpha/2}$,
 which is  the transition density of 
 the rotationally symmetric $\alpha$-stable process $Y$,
  by the following Duhamel's formula:
 \begin{equation} \label{e:5.1}
  q^b (t, x, y)
 =p (t, x, y)+
 \int^t_0\int_{\R^d}q^b (s, x, z)\, b(z)\cdot\nabla_z p (t-s, z, y)
 d z  ds.
 \end{equation}
 Applying the above formula recursively, one expects that
     \begin{equation} \label{e:5.2}
     q^b(t, x, y):= \sum_{k = 0}^{\infty} q^b_k (t, x, y)
   \end{equation}
   is a fundamental solution for $\LL^b$,
   where
  $q^b_0(t, x, y) := p(t, x, y)$ and for $k\geq 1$,
 \begin{equation} \nonumber
 q^b_k (t, x, y):= \int_0^t \int_{\R^d} q^b_{k-1} (s, x, z)\, b(z) \cdot
 \nabla_z p(t-s, z, y) d z.
 \end{equation}
It is shown in \cite{BJ} that the series in \eqref{e:5.2} converges absolutely and, for every $T>0$, such defined
$q^b(t, x, y)$ is a conservative transition density function and  is comparable to $p(t, x, y)$ on $(0, T]\times \R^d \times \R^d$.
Recall that $p(t, x, y)$ has two-sided estimate \eqref{e:2.5}.
In \cite{CLW}, Chen and Wang showed that the Markov process $X_t$ having $q^b(t, x, y)$ as its transition density function
is the unique solution to the martingale problem $(\LL^b, C^2_b (\R^d)$; moreover, it is the unique weak solution to the
following stochastic differential equation:
$$ dX_t =dY_t +b(X_t) dt,
$$
where $Y_t$ is the rotationally symmetric $\alpha$-stable process on $\R^d$.
Dirichlet heat kernel estimate for  
${\cal L}^b$ in a bounded $C^{1,1}$ open set
has been obtained in \cite{CKS1}.
In \cite{KS1, KS2}, Kim and Song extended results in \cite{BJ, CKS2} to $\Delta^{\alpha/2} + \mu \cdot \nabla$,
where $\mu=(\mu_1, \dots, \mu_d)$ are signed measures in suitable Kato class.
These work can be regarded as heat kernels for fractional Laplacian under gradient perturbation.
Heat kernel estimates for relativistic stable processes and for mixed Brownian motions and stable processes 
 with drifts
have recently been studied in \cite{CLW2} and \cite{CH}, respectively.
See \cite{CD} for drift perturbation of  subordinate Brownian motion of pure jump type and its heat kernel estimate.
 While  in \cite{Xie-Zhang},
 Xie and Zhang considered the critical operator $\sL^b:=a\Delta^{1/2}+b\cdot \nabla$, where 
for some $0<c_0<c_1$, $a:\R^d\to[c_0,c_1]$ and $b:\R^d\to\R^d$ are two H\"older continuous functions. 
They established two-sided estimates for the heat kernel of $\sL^b$ by using Levi's method as described in Subsection 3.1.

In the same spirit, Wang and Zhang in \cite{Wang-Zhang} considered more general fractional diffusion operators  
over a complete Riemannian manifold perturbed by a time-dependent gradient term, and showed two-sided estimates and gradient estimate of the heat kernel.
More precisely, let $M$ be a $d$-dimensional connected complete Riemannian manifold with Riemannian distance $\rho$. Let $\Delta^M$ be the Laplace-Beltrami operator. Suppose that the heat kernel $p(t,x,y)$ of $\Delta^M$ with respect to the Riemannian volume $dx$ 
exists and has the following two-sided estimates:
\begin{align}\label{ET1}
c_1t^{-d/2}e^{-c_2\rho(x,y)^2/t}\leq p(t,x,y)\leq c_3t^{-d/2}e^{-c_4\rho(x,y)^2/t},\ \ t>0, x,y\in M,
\end{align}
and gradient estimate
\begin{align}\label{ET2}
|\nabla_xp(t,x,y)|\leq c_5t^{-(d+1)/2}e^{-c_4\rho(x,y)^2/t},
\end{align}
where $\nabla_x$ denotes the covariant derivative.
Let $P_t$ be the corresponding semigroup, that is,
$$
P_t f(x):=\int_M p(t,x,y)f(y)dx,\ \ f\in C_b(M).
$$
For $0<\alpha< 2$, consider the $(\alpha/2)$-stable subordination of $P_t$
$$
P_t^{(\alpha)}:= \int_0^\infty P_s\,\mu_t^{(\alpha/2)}(ds),\quad  t\ge 0,
$$
where $\mu_t^{(\alpha/2)}$ is a probability measure on $[0,\infty)$ with Laplace transform
$$
\int_0^\infty e^{-\lambda s} \mu_t^{(\alpha/2)}(ds)=e^{-t\lambda^{\alpha/2}}, \quad \lambda\geq 0.
$$ 
Then $P_t^{(\alpha)}$ is  a $C_0$-contraction  semigroup on $C_b(M)$.  Let $\LL^{(\alpha)}$ be  the infinitesimal generator of $P_t^{(\alpha)}$.
In \cite{Wang-Zhang}, Wang and Zhang considered the following operator
$$
 \LL_{b,c}^{(\alpha)}f(t,x) := \LL^{(\alpha)}f(x) +\< b(t,x), \nabla_xf(x)\> + c(t,x)f(x),\ f\in C^2_b(M),
$$
where $b:\R_+\times M\to TM$ and $c:\R_+\times M\to\R$ are measurable. For $\alpha\in(0,2)$, one says that a measurable function 
$f:\R_+\times M\to \R$ belongs to Kato's class $\mK_\alpha$ if
$$
\lim_{\varepsilon\to 0}\sup_{(t,x)\in [0,\infty)\times M} \varepsilon^{1/\alpha} \int^\varepsilon_0\!\!\!\int_M 
 \frac{s^{1-1/\alpha}(\varepsilon-s)^{-1/\alpha}|f(t\pm s, y)|}{(s^{1/\alpha}+\rho(x,y))^{d+\alpha}}dyds =0.
$$
Notice that when $\alpha=2$ and $f$ is time-independent, $\mK_\alpha$ is the same as in \eqref{ET3}.

The following result is shown
in \cite{Wang-Zhang}.
 \begin{thm}
Assume \eqref{ET1}, \eqref{ET2} and $\alpha\in (1,2).$
 If $|b|,c\in \mK_\alpha$, then there is a unique continuous function $p_{b,c}^{(\alpha)}(t,x;s,y)$ having the following properties:
\begin{enumerate}

\item[\rm (i)] (Two-sided estimates)  There is a constant $c_1>0$ such that for all $t-s\in (0,1], x,y\in M$,
$$
c_1^{-1}\frac{t-s}{((t-s)^{1/\alpha}+\rho(x,y))^{d+\alpha}}\leq p_{b,c}^{(\alpha)}(t,x;s,y)\leq c_1 \frac{t-s}{((t-s)^{1/\alpha}+\rho(x,y))^{d+\alpha}}.
$$
 \item[\rm (ii)] (Gradient estimate) There is a constant $c_2>0$ such that for all $t-s\in (0,1], x,y\in M$,
$$
|\nabla_x p_{b,c}^{(\alpha)}(t,x; s,y)|\leq c_2\frac{(t-s)^{1-1/\alpha}}{((t-s)^{1/\alpha}+\rho(x,y))^{d+\alpha}}.
$$ 

\item[\rm (iii)] (C-K equation) For any $0\leq s<r<t$ and $x,y\in M$,
$$
p_{b,c}^{(\alpha)}(t,x;s,y)= \int_M p_{b,c}^{(\alpha)}(t,x;r,z)p_{b,c}^{(\alpha)}(r,z;s,y)dz.
$$
\item [\rm (iv)] (Generator) If $b\in C([0,\infty); L_{loc}^1(M,dx; TM))$ and $c\in C([0,\infty); L_{loc}^1(M,dx; \R))$, 
then for any $\varphi,\psi\in C_0^2(M),$
$$
\lim_{t\downarrow s} \frac 1 {t-s} \int_M \psi(P_{t,s}^{b,c}  \varphi -\varphi)dx=\int_M\psi \LL_{b,c}^{(\alpha)}(s,\cdot)\varphi dx,\ \ s\ge 0,
$$  
where $P_{t,s}^{b,c}\varphi:= \int_M p_{b,c}^{(\alpha)}(t,\cdot; s,y)\varphi(y) dy.$
\end{enumerate} 
\end{thm}
 
\medskip

The above results indicate that, under suitable Kato class condition, heat kernel estimates are stable under gradient perturbation.

 In \cite{CW}, Chen and   Wang   studied
 heat kernels for fractional Laplacian under  non-local perturbation of high intensity;
 that is, heat kernels for
\begin{equation}\label{e:5.3}
\LL^\kappa f(x)=\Delta^{\alpha/2} f(x) + \SS^\kappa f(x),
\quad f\in C^2_b(\R^d),
\end{equation}
where
 \begin{equation}\label{e:5.4}
\SS^\kappa f(x):={\cal A}(d, -\beta)
\int_{\R^d} \left( f(x+z)-f(x)- \nabla f(x) \cdot
z \1_{\{|z|\leq 1\}}  \right) \frac{\kappa(x, z)}{|z|^{d+\beta}}dz
\end{equation}
for some $0<\beta <\alpha <2$ and   a  real-valued bounded function $\kappa(x, z)$
on $\R^d\times \R^d$ satisfying
$$
\kappa(x, z)=\kappa(x, -z) \qquad \hbox{for every } x, z\in \R^d.
$$
 Uniqueness and existence of fundamental solution $q^\kappa (t, x, y)$ is established in \cite{CW}.
 The approach is also a perturbation argument by viewing   $\LL^\kappa =\Delta^{\alpha/2}+\SS^\kappa $ as  a lower order perturbation of
 $\LL^0=\Delta^{\alpha/2}$ by $\SS^\kappa $.
 So  heuristically,  the fundamental solution
(or heat kernel) $q^\kappa (t, x, y)$ of $\LL^\kappa $ should satisfy the following
Duhamel's formula:
\begin{equation}\label{e:5.5}
q^\kappa (t, x, y)=p (t, x, y)+\int_0^t \int_{\R^d} q^\kappa  (t-s, x, z)
\SS^\kappa _z p  (s, z, y)dz ds
\end{equation}
for $t>0$ and $x, y\in \R^d$.
Here the notation $S^\kappa _z p  (s, z, y)$ means that the non-local operator $S^\kappa $
is applied to the function $z\mapsto p(s, z, y)$.
Similar notation will also be used for other operators, for example, $\Delta^{\alpha/2}_z$.
Applying \eqref{e:5.5} recursively, it is reasonable to
conjecture that $\sum_{n=0}^\infty q^\kappa _n(t, x, y)$, if convergent,
 is a solution to \eqref{e:5.5}, where $q^\kappa _0(t, x, y):=p  (t, x, y)$
 and
\begin{equation}\label{e:qn}
q^\kappa _n (t, x, y):=\int_0^t \int_{\R^d} q^\kappa _{n-1} (t-s, x, z)\SS^\kappa _z p(s, z, y) dz ds \quad \hbox{for } n\geq 1.
\end{equation}
The hard part is the estimates on $\SS^\kappa _z p  (s, z, y)$ and on each $q^\kappa _n(t, x, y)$.
 In contrast to the gradient perturbation case, 
 the fundamental solution
 to the non-local perturbation \eqref{e:5.3} does not need to be positive, and when the kernel is positive,
 it does not need to be comparable to $p(t, x, y)$.
 One can rewrite $\LL^\kappa $ of \eqref{e:5.3} as follows:
 $$ \LL^\kappa f (x)=\int_{\R^d} \left( f(x+z)-f(x)-\langle\nabla f(x), z \>
\1_{\{|z|\leq 1 \}} \right) j^\kappa (x, z)dz,
$$
where
\begin{equation}\label{e:5.6}
j^\kappa (x, z)= \frac{{\cal A}(d, -\alpha)}{|z|^{d+\alpha}}
\left( 1+ \frac{{\cal A}(d, -\beta)}{{\cal A}(d, -\alpha)}\,
\kappa(x, z)\, |z|^{\alpha-\beta} \right) .
\end{equation}
 It is shown in \cite{CW} that the fundamental solution kernel $q^\kappa \geq 0$ if  $j^\kappa (x, z)\geq 0$;
 that is, if
 \begin{equation}\label{e:5.7}
\kappa(x, z)\geq -\frac{{\cal A} (d, -\alpha)}{{\cal A} (d, -\beta)} \,
|z|^{\beta -\alpha} \quad
\hbox{for  a.e. } z\in \R^d.
\end{equation}
When $\kappa(x, z)$ is continuous in $x$, the above condition is also necessary for the non-negativity of $q^\kappa (t, x, y)$.
Under condition \eqref{e:5.7}, various sharp heat kernel estimates have been obtained in \cite{CW}.
In particular, it is shown in \cite{CW} that if there are constants $0<c_1\leq c_2$ so that
$$
\frac{c_1}{|z|^{d+\alpha} }\leq j^\kappa (x, z) \leq \frac{c_2}{|z|^{d+\alpha} } \quad \hbox{for } x, z\in \R^d,
$$
then for every $T>0$, $q^\kappa (t, x, y) \asymp p(t, x, y)$ on $(0,T]\times \R^d\times \R^d$.
Dirichlet heat kernel estimates for $\LL^\kappa $ of \eqref{e:5.3} has recently been studied in Chen and Yang \cite{CY}.

In a subsequent work \cite{Wang},  Wang studied fundamental solution for $\Delta +\SS^\kappa $ and its two-sided heat kernel estimates.
In \cite{CKS2}, Chen, Kim and Song established stability of heat kernel estimates under (local and non-local) Feynman-Kac
transforms for a class of jump processes; see also C. Wang \cite{W} on a related work. 
Very recently, stability of heat kernel estimates for diffusions with jumps (both symmetric and non-symmetric) under Feynman-Kac transform has been studied in Chen and  Wang \cite{CWa}.
 On the other hand,
 by employing the strategy and road map from Chen and Zhang \cite{CZ} as outlined in
Section 3 of this paper, 
Kim, Song and Vondracek \cite{KSV} has extended Theorem \ref{T:3.1} to more general 
non-local operator
 $\sL$ of \eqref{e:3.1} with $\frac{1}{|z|^{d+\alpha}}$ being replaced by the density of
 L\'evy measure of certain
 subordinate Brownian motions. In a recent work \cite{CCW}, X. Chen, Z.-Q. Chen 
and J. Wang have used Levi's freezing coefficient method to obtain
upper and lower bound estimates
for heat kernels of the following type of non-local operators of variable order:
$$
\LL f(x):=\int_{\R^d}\Big(f(x+z)-f(x)-\nabla f(x)\cdot z \1_{\{|z|\le
1\}}\Big) \frac{\kappa(x,z)}{|z|^{d+\alpha(x)}}\,dz, \quad f\in C^2_c (\R^d), 
$$
where $\alpha (x)$ is a H\"older continuous function on $\R^d$ such that 
$$
0<\alpha_1\leq \alpha (x) \leq \alpha_2<2 \quad \hbox{for all } x\in \R^d,
$$
 and $\kappa (x, z)$ satisfies conditions \eqref{e:2.6}-\eqref{e:2.7}.

In this survey, we mainly concentrate on the quantitive estimates of the heat kernels of non-symmetric nonlocal operators.
For  derivative formula of the heat kernel associated with stochastic differential equations with jumps, 
we refer the interested reader to 
\cite{Zhang0, Wang-Xie-Zhang, Wang-Xu-Zhang}.
For other results on the existence and smoothness of heat kernels or fundamental solutions 
for non-symmetric jump processes or non-local operators under H\"ormander's type conditions,  see 
\cite{PZ, Ku00} 
for the studies of linear Ornstein-Uhlenbeck processes with jumps, and 
\cite{Zhang1, Zhang2, Zhang3}
and the references therein for the studies of general stochastic differential equations with jumps. We will not survey these results since the arguments in   the above references are mainly based on the Malliavin calculus
and thus belong to another topic. 

\bigskip

{\bf Acknowledgement}. We thank the referee for helpful comments, in particular for pointing out a gap in the proof of (2.30) in \cite{CZ}.

\vskip 0.3truein

\noindent {\bf Zhen-Qing Chen}

\noindent Department of Mathematics, University of Washington,
Seattle, WA 98195, USA

\noindent E-mail: zqchen@uw.edu

\bigskip

\noindent {\bf Xicheng Zhang}

\noindent School of Mathematics and Statistics, Wuhan University, Hubei 430072, P. R. China

\noindent E-mail: XichengZhang@gmail.com

\end{document}